\documentclass[11pt]{article}

\usepackage{graphicx}

\usepackage[utf8]{inputenc}
\usepackage{amsthm}

\begin{document}

\begin{center}

{\Large \noindent \bf
Connected and Autonomous Vehicle Scheduling Problems: Some Models and Algorithms.}

\vspace{7mm}

 {\bf  Evgeny R. Gafarov}
\medskip

{\small \it V.A. Trapeznikov Institute of Control Sciences of the Russian Academy of Sciences,\\
Profsoyuznaya st. 65, 117997 Moscow, Russia, \\
 email:  axel73@mail.ru }

\medskip

{\bf Frank Werner}
\medskip

{\small \it Fakult\"at f\"ur Mathematik, Otto-von-Guericke-Universit\"at Magdeburg,\\
PSF 4120, 39016 Magdeburg, Germany, \\ email: frank.werner@mathematik.uni-magdeburg.de}

\vspace{5mm}
{\it April 04, 2023 }

\vspace{3mm}

\end{center}

\abstract{%
In this paper, we consider scheduling problems that arise in connected and autonomous vehicle systems. For four variants of such problems, mathematical models and solution algorithms are presented. In particular, three polynomial algorithms and a branch and bound algorithms are developed.
}



\bigskip

\newpage

\section{Introduction}

A car-to-car communication has been presented as a possible solution to many challenges encountered in this field. Many solutions have been presented involving the modeling of the entirety of an autonomous driving system as a muli-agent system, where vehicles interact to enable autonomous functionality such as emergency braking and traffic jam avoidance.
Vehicle systems are developing towards fully connected and fully autonomous systems.
Vehicular communication technologies have been considered e.g. in \cite{Bazzal,Sahin}.

These vehicle systems can coordinate the vehicles in order to speed up the traffic and avoid traffic jams.  Vehicles can be coordinated by a centralized scheduler residing in the network (e.g., a base station in case of cellular systems) or a distributed scheduler, where the resources are autonomously selected by the vehicles

In \cite{Qian}, the authors proposed to optimize the departure times, travel routes, and longitudinal trajectories of connected and autonomous vehicles (CAVs) as well as signal timings at intersections to achieve a stable traffic state, where no vehicles need to stop before entering any intersection and no queue spillover occurs at any intersection. The departure times, travel routes and signal timings are optimized in a central controller, while the vehicle trajectories can be optimized by distributed roadside processors, which together form a hierarchical traffic management scheme.

In \cite{Atagoziev}, the authors considered the coordination of lane changes of autonomous vehicles on a two-lane road segment before reaching a given critical position. An algorithm is presented that performs a lane change of a single vehicle in the shortest possible time. This algorithm is then iteratively applied in order to handle all lane changes required on the considered road segment while guaranteeing traffic safety.

In \cite{Ma}, the scheduling problem of a CAV
crossing the intersection was considered to optimize the intersection efficiency. In addition, a solution algorithm was presented. 

In \cite{Younes}, the  time phases of the traffic light scheduling problem were considered with the goal of increasing the traffic fluency by decreasing the waiting time of the traveling vehicles 
at the signalized road intersections.

In this paper, we consider four scheduling problems that arise in connection with CAVs. The reminder of this paper is as follows. In each of the Sections 2 - 5, we consider one of these problems. The problem with a road of two lanes and a barrier on the lane is considered in Section 2. Section 3 deals with the case of a turn to a main road. Section 4 considers the case of a road with three lanes and a barrier on the middle lane. Section 5 deals with a crossroad having dividing lanes. 
For each of these cases, an appropriate scheduling problem is formulated and a solution algorithm is given. Finally, Section 6 gives a few concluding remarks. 


\section{A road with two lanes and a barrier on a lane} 

\begin{figure}
\begin{center}
  \includegraphics[width=0.7\textwidth]{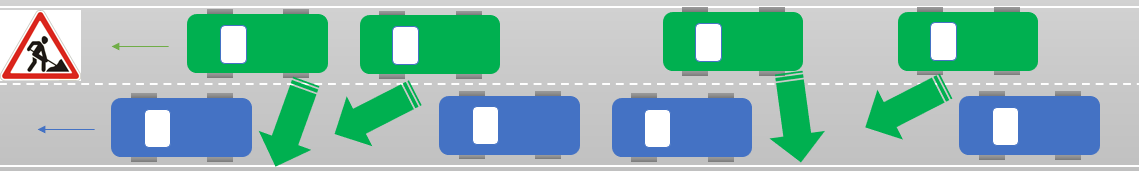}
 \caption{A road with two lanes and a barrier on a lane}\label{fig1}
\end{center}
\end{figure}

In this section, we consider a road with two lanes, where two sets $N_1$ and $N_2$ of CAVs are given. The CAVs from the set $N_1$ go on lane 1, and the CAVs from the set $N_2$ go on lane 2.  Both lanes have the same direction. On lane 2, there is a barrier and the CAVs from the set $N_2$ have to move to lane 1, see Fig. \ref{fig1}.

We have to find a sequence of passing the barrier by the CAVs from the sets $N_1$ and $N_2$ in order to minimize a given objective function, e.g. the total passing time.

We assume that
\begin{itemize}
\item a maximal feasible speed of the CAVs is given. The CAVs either go with the maximal feasible speed or brake in order to let another CAV change the lane;
\item an acceleration is not taken into account;
\item te time needed to change the lane is not taken into account, i.e., it is equal to zero;
\item the CAVs have the same length;
\item the safe distance between two CAVs is the same for all vehicles.
\end{itemize}

The same problem arises e.g. on railway sections and in automated warehouses of logistics companies with autonomous robot-transporters. 
This simplified problem can be formulated as a single machine scheduling problem as follows.

Given a set $N = N_1 \bigcup N_2$ of $n$ jobs that have to be processed on a single machine from time $0$ on. For each job $j$, a processing time $p_j = p$, a release date $r_j>0$, a due date $d_j$ and a weight $w_j$ are given.
The processing time $p$ can be computed from the maximal feasible speed and the length of a CAV. The value $r_j$ corresponds to the position of the CAV $j$ on the road.

A schedule is uniquely determined by a permutation $\pi$ of the CAVs of the set $N$. 
Let $C_{j}(\pi)$ be the completion time of job $j$ in the schedule $\pi$. A precedence relation can be defined, i.e., for the jobs from the set $N_1$, we have $j^1_1 \rightarrow j^1_2 \rightarrow \dots \rightarrow j^1_{n_1}$, where $n_1 = |N_1|$ and $j\rightarrow i$ means that the processing of job $j$ precedes the processing of job $i$. Thus, there is a chain of jobs on lane 1. Analogously, a chain of jobs can be defined for the set $N_2$.

For the single machine scheduling problem of minimizing total completion time, the goal is to find an optimal schedule $\pi^*$ that minimizes the total completion time, i.e., 
\begin{equation}
\sum C_j = \sum_{j\in N} C_{j}.
\end{equation}

Here the completion time of a job is equal to the time when the car passes the barrier. We denote this problem by $1|2\ chains, p_{j}=p, r_j|\sum C_j$ according to the traditional three-field notation $\alpha|\beta|\gamma$ for scheduling problems proposed by Graham et al.~\cite{Graham79}, where $\alpha$ describes the machine environment, $\beta$ gives the job characteristics and further constraints, and $\gamma$ describes the objective function.

Let $$T_j(\pi) = \max\{0,C_{j}(\pi)-d_j\}$$ be the tardiness of job $j$ in the schedule $\pi$. In addition, one can consider also the following  objective functions: 
\begin{eqnarray*}
\sum w_jC_j & = & \sum w_jC_j(\pi) \hbox{ -- total weighted completion time,}  \\
\sum T_j & = &   \sum T_j(\pi) \hbox{ -- total tardiness,}  \\  
\sum w_j T_j & = &  \sum w_jT_j(\pi) \hbox{ -- total weighted tardiness.}   
\end{eqnarray*}
It is known that the problems $1|chains, p_{j}=p, r_j|\sum w_jC_j$ and $1|chains, p_{j}=p, r_j|\sum w_jT_j$ with an arbitrary number of chains are  NP-hard \cite{Lenstra}. This has been proven by a reduction from the 3-Partition Problem.

In \cite{Baptiste2000}, Baptiste presented polynomial time dynamic programming algorithms to solve the problems $1|p_j = p, r_j|\sum T_j$ and $1|p_j = p, r_j|\sum w_jU_j$.  

In an optimal schedule for the problem $1|2\ chains, p_{j}=p, r_j|\sum C_j$, the jobs are processed in non-decreasing order of the values $r_j$. This can be easily proven  by  contradiction. For an illustration of the concepts introduced above, we consider the following small example.

{\bf Example.}  Let $N_1 = \{1,2\},\ N_2=\{3,4\}$. Moreover, the values $p = 2,\ r_1 = 0,\ r_2 = 3,\ r_3 = 1,\ r_4=4$ and $d_1=d_2=0,\ d_3=3,\ d_4=6$ are given. For the chosen job sequence $\pi = (1,3,2,4)$, we obtain $S_1(\pi)=0,\ S_3(\pi)=2,\ S_2(\pi)=4,\ S_4(\pi) = 6$,\ $C_1(\pi)=2,\ C_3(\pi)=4,\ C_2(\pi)=6,\ C_4(\pi) = 8$. $\sum_{j=1}^4 C_j(\pi) = 20$ and $\sum_{j=1}^4 T_j(\pi) = 1+2=3$. For the job sequence $\pi' = (3,4,1,2)$, we get  $\sum_{j=1}^4 T_j(\pi') = 0$.

We note that there exists a set $\Theta$ of possible completion times of all jobs, where $|\Theta| \leq n^2$, since:
\begin{itemize}
\item without loss of generality, we consider only active schedules, where no job can be processed earlier without loss of feasibility;
\item there are no more than $n$ different values $r_j$;
\item all processing times are equal to $p$ and thus, for any job $j \in N$, its completion time is equal to $r_i+lp,\ i \in N, l\leq n$.
\end{itemize}

\bigskip
{\bf Theorem 1.} The problems $1|2\ chains, p_{j}=p, r_j|f,\ f\in\{\sum w_jC_j, \sum w_jT_j\}$ can be solved in $O(n^5)$ time by a dynamic program.  \\

A sketch of the proof is as follows. In the dynamic program (DP), we consider the jobs $i_1,i_2,\dots,i_{n_2} \in N_2$ one by one, where  $i_1 \rightarrow i_2 \rightarrow \dots \rightarrow i_{n_2}$. Thus, at each stage $k$ of the dynamic program, we consider a single job $i_k,\ k=1,2,\dots,n_2$. Moreover, at each stage we consider all states $(f, C_{max}, pos)$ stored at the previous stage. In addition, for each state, we store the best partial solution (sequence of jobs). The meaning of the above triplet is as follows. Here $pos \in \{0,1,2,\dots,n_1\}$ describes the position of a job, and it means that job $i_{k-1}$ is processed between the jobs $j_{pos} \in N_1$ and $j_{pos+1}\in N_1,\ 0<pos<n_1$, and $C_{max} = C_{i_{k-1}}\in \Theta$ denotes the completion time of job $i_{k-1}$ in the corresponding partial solution. Finally, $f$ is the value of the considered objective function that corresponds to the partial solution. For each job $i_k$ and a state $(f, C_{max}, pos)$, we compute new states $(f', C'_{max}, pos')$, where $pos'\geq pos$, and $C'_{max}$ is the completion time of job $i_{k}$ in a new partial solution, where job $i_k$ is scheduled after job $j_{pos'} \in N_1$.
If at any stage, there are two states $(f', C'_{max}, pos')$ and $(f'', C''_{max} pos')$ with $f'\leq f''$ and $C'_{max} \leq C''_{max}$, we only keep the state  $(f', C'_{max}, pos')$.
After the last stage, we have to select the best found complete solution among all states generated.

A pseudo-code of Algorithm DP is presented below. \\

{\bf Algorithm DP} 

\begin{list}{}{}
\item[1.]  $StatesSet = \{(0,0,0)\}$;
\item[2.] FOR EACH $i_k \in N_2$ DO
\begin{list}{}{}
\item[2.1] $NewStatesSet = \{\}$;
\item[2.2] FOR EACH $(f, C_{max}, pos)\in StatesSet$ DO
           \begin{list}{}{}
              \item[2.2.1] Let $PositionsList = \{pos,pos+1,\dots, n_1\}$;
              \item[2.2.2] FOR EACH $pos' \in PositionsList$ DO
              \begin{list}{}{}
                  \item[2.2.2.1] Calculate $f'$ for the resulting partial solution, if job $i_k$ is processed after $j_{pos'}$, according to the partial solution corresponding to state $(f,  C_{max}, pos)$;
                  \item[2.2.2.2] Add  $(f',  C'_{max}, pos')$ to $NewStatesSet$. If in $NewStatesSet$, there is a state $(f'',  C''_{max}, pos')$ with $f'\leq f''$ and $C'_{max} \leq C''_{max}$, then exclude the state  $(f'',  C''_{max}, pos')$ from  $NewStatesSet$;
              \end{list}
           \end{list}
\item[2.3] $StatesSet :=NewStatesSet$;           
\end{list} 
\item[3.] Select the best found complete solution among all states generated.
\end{list}

\section{Turn to a main road}

\begin{figure}
\begin{center}
  \includegraphics[width=0.7\textwidth]{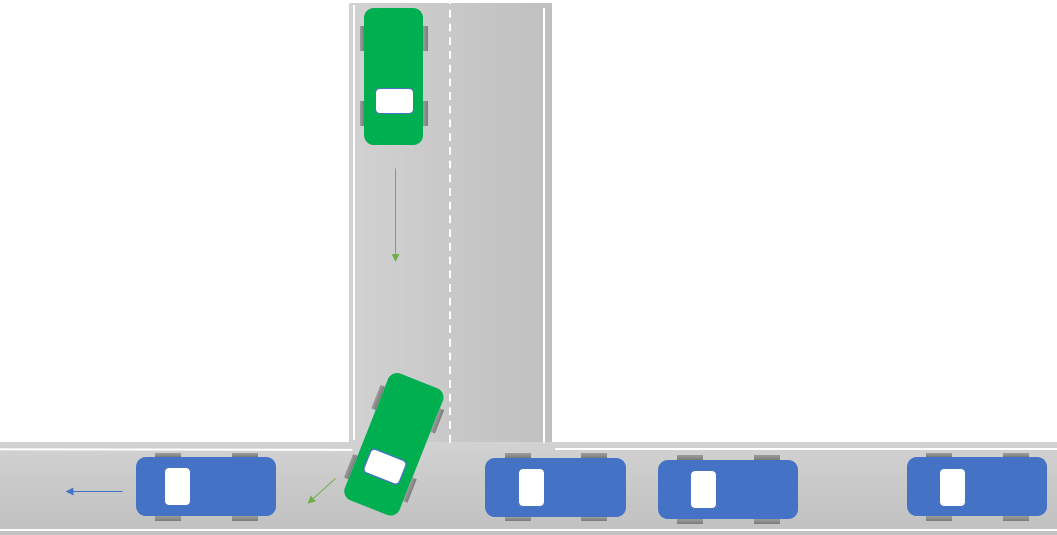}
 \caption{Turn to a main road}\label{fig2}
\end{center}
\end{figure}

There is a set $N_1$ of CAVs  going along a main road and a set $N_2$ of CAVs  turning into the main road from a side road, see Fig. \ref{fig2}. 
In contrast to  the problems $1|2\ chains,p_{j}=p|\gamma$, we have now $p_j=p^1,\ j \in N_1$ and $p_j=p^2,\ j \in N_2$. We denote this problem by $1|2\ chains, p_{j}\in {p^1,p^2}, r_j|\gamma$. This problem can be solved by the same DP.

\section{A road with three lanes and a barrier on the middle lane}

\begin{figure}
\begin{center}
  \includegraphics[width=0.7\textwidth]{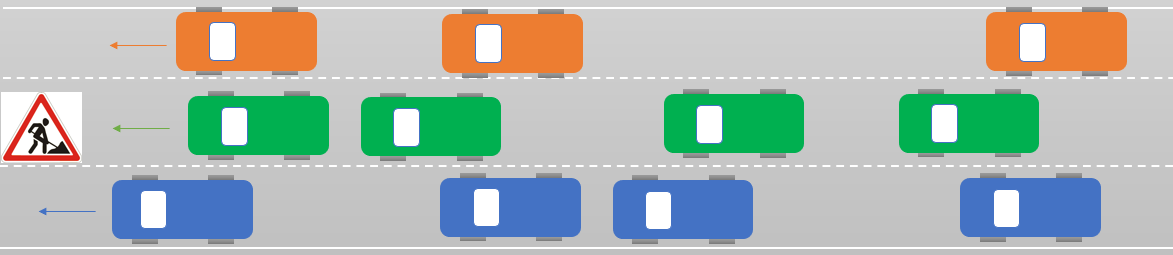}
 \caption{A road with three lanes and a barrier on the middle lane}\label{fig3}
\end{center}
\end{figure}

In addition to the problems $1|2\ chains,p_{j}=p|\gamma$, there are an additional lane $3$  and a subset $N_3$ of jobs, see Fig. \ref{fig3}. The jobs of the set $N_1$ should be processed on the machine $1$, and the jobs of th set $N_3$  should be processed on the machine $3$. The jobs of the set $N_2$ can be processed on any of these two machines. 
Precedence relations among the jobs of the set $N_3$ can be defined as a chain of jobs.

We denote this problem by $P2|dedicated,3\ chains, p_{j}=p, r_j|\gamma$.
This problem can be solved by a modified DP, where we consider the positions $pos$  between the jobs of the set $N_1$ and between the jobs of the set $N_3$.

\section{A crossroad with dividing lines}

\begin{figure}
\begin{center}
  \includegraphics[width=0.5\textwidth]{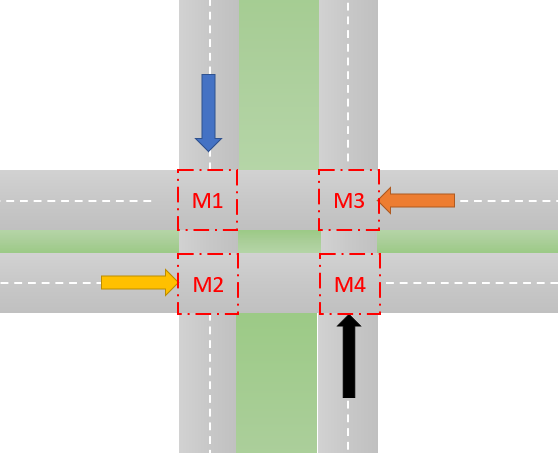}
 \caption{A crossroad with dividing lines}\label{fig4}
\end{center}
\end{figure}

In this section, we consider a crossroad with dividing lines and four sets  $N_1, N_2, N_3, N_4$ of CAVs. They share four sectors of the crossroad denoted by $M_1, M_2, M_3, M_4$, see Fig. \ref{fig4}.  We have to find an optimal sequence of passing these sectors.
 
We can formulate the following job shop scheduling problem with four machines. There are four sets $N_1, N_2, N_3, N_4$ of jobs  and four machines corresponding to the sectors $M_1, M_2, M_3, M_4$. Each job $j$ consists of two operations. For each job $j\in N_1$, its first operation has to be processed on machine $M_1$ and its second one has to be processed on machine $M_2$.
For each job $j\in N_2$, its first operation has to be processed on machine $M_2$ and its second one has to processed on machine $M_4$.
For each job $j\in N_3$, its first operation has to be processed on machine $M_3$ and its second one has to be processed on machine $M_1$.
For each job $j\in N_4$, its first operation has to be processed on machine $M_4$ and its second one has to be processed on machine $M_3$. The processing times of the operations are equal to $p$. 
Precedence relations can be given as chains of jobs.

If the lengths of the dividing lines are equal to $0$, then the second operation of a job $j$ should be processed immediately after the first one. Otherwise for each of the sets $N_1, N_2, N_3, N_4$, there are four buffers of limited capacities, namely $b_1, b_2, b_3, b_4$ jobs for the corresponding machine. At any moment, for the set $N_1$, there can be up to $b_1$ jobs for which the first operation is completed and the second one is not yet started. 
We denote these problems by $J4|4\ chains, p_{j}=p, r_j|f,$\\ $f\in\{C_{max}, \sum w_jC_j, \sum w_jT_j\}$, where $C_{max}$ is the makespan.

The problems $J4|4\ chains, p_{j}=p, r_j|f,\ f\in\{C_{max}, \sum w_jC_j, \sum w_jT_j\}$ can be solved by a branch-and-bound (B\&B) algorithm. The search (rooted) tree is constructed by the following branching rule. For any node of the tree, we consider the following 8 possible branches:
\begin{itemize}
\item Schedule the first unscheduled possible operation for a job $j\in N_1$ on machine $M_1$ at the earliest possible starting time. If there is no such an operation, skip this branch.
\item Schedule the first unscheduled possible operation for a job $j\in N_3$ on machine $M_1$ at the earliest possible starting time. 
\item Schedule the first unscheduled possible operation for a job $j\in N_1$ on machine $M_2$ at the earliest possible starting time. 
\item Schedule the first unscheduled possible operation for a job $j\in N_2$  on machine $M_2$ at the earliest possible starting time. 
\item Schedule the first unscheduled possible operation for a job $j\in N_3$  on machine $M_3$ at the earliest possible starting time. 
\item Schedule the first unscheduled possible operation for a job $j\in N_4$  on machine $M_3$ at the earliest possible starting time. 
\item Schedule the first unscheduled possible operation for a job $j\in N_2$  on machine $M_4$ at the earliest possible starting time. 
\item Schedule the first unscheduled possible operation for a job $j\in N_4$  on machine $M_4$ at the earliest possible starting time. 
\end{itemize}

Thus, there are up to $2^3 = 8$ branches for each node to be considered. Since there are $2n$ operations, where $n=|N_1  \bigcup N_2  \bigcup N_3  \bigcup N_4|$, there are no more than $2n$ levels in the search tree.  Thus, we have no more than $(2^3)^{2n} = 2^{6n}$ nodes to be considered. If some of the values $b_1, b_2,b_3,b_4$ are equal to $0$, we have less nodes, e.g., if each of them is equal to $0$, then we have only $2^{3n}$ nodes.

Moreover, we can use the following trivial lower and upper bounds for the problem $J4|4\ chains, p_{j}=p, r_j|C_{max}$.  \\

{\bf Upper bound.} To construct a feasible solution, we use a list scheduling algorithm. In this algorithm, we consider the unscheduled operations one-by-one according to a non-decreasing order of the release dates of the corresponding jobs. We schedule the next unscheduled operation at the earliest possible starting time according to the current partial schedule. To order the set of jobs, we need $O(n\log n)$ operations. In addition, we need $O(n)$ operations to construct a feasible solution.  \\

{\bf Lower bound.} Consider a set of unscheduled operations $N'$. For each of them, we calculate the earliest possible starting time according to the current partial schedule without taking into account the other unscheduled operations. In such a way, we get a schedule $\pi$ that can be infeasible. Let $C_{M_1}(\pi)$ be the makespan (i.e., the maximal completion time of an operation assigned to the machine) for machine $M_1$, $IT_{M_1}(\pi)$ be the idle time of machine $M_1$ between the operations of the set $N'$, and $OT_{M_1}(\pi)$ be the total overlap time, where more than one operation is processed at the same time. Moreover, let 
$$C'_{M_1}(\pi) = C_{M_1}(\pi)+\max\{0,OT_{M_1}(\pi)-IT_{M_1}(\pi)\}. $$ Then 
$$LB1 = \max\{C'_{M_1}(\pi), C'_{M_2}(\pi), C'_{M_3}(\pi), C'_{M_4}(\pi)\}$$ is a lower bound. It is easy to check that we need $O(n)$ operations to calculate this bound.  \\

If we use Upper bound and Lower bound, then the B\&B algorithm requires $O(n2^{6n})$ operations. 

\section{Concluding Remarks}

In this note, four models of scheduling problems for CAVs have been given. Three of them can be solved by a dynamic programming algorithm in polynomial time. For the fourth problem, a B\&B algorithm has been presented.

The following questions can be considered in the future:.
\begin{itemize}
\item Are the problems $J4|4\ chains, p_{j}=p, r_j|f,$  $f\in\{C_{max}, \sum w_jC_j, \sum w_jT_j\}$ NP-hard or can be solved in polynomial time?
\item Are there problems with CAVs having equal processing times and a fixed number of chains of jobs which is an NP-hard problem?
\end{itemize}

\end{document}